\documentclass[12pt,twoside]{article}
\usepackage[english]{babel}
\usepackage{amsmath}
\usepackage{times,amssymb,amscd}
\usepackage[utf8]{inputenc}
\usepackage[pass,legalpaper]{geometry}
\usepackage{fancyhdr}
\usepackage{multicol}
\usepackage{enumitem}
\usepackage{graphicx}
\usepackage{pgf}
\usepackage{textcomp}
\newtheorem{thm}{Theorem}[section]

\newtheorem{defn}[thm]{Definition}

\numberwithin{equation}{section}

\newcommand{\bG}{\mathbf{G}}

\newcommand{\ba}{\mathbf{a}}

\newcommand{\bh}{\mathbf{h}}

\newcommand{\Bc}{\boldsymbol{c}}

\newcommand{\Bu}{\boldsymbol{u}}

\newcommand{\cD}{\mathcal{D}}

\newcommand{\cT}{\mathcal{T}}
\newcommand{\cC}{\mathcal{C}}

\newcommand{\cB}{\mathcal{B}}

\begin{document}
\pagestyle{myheadings}
\markboth{\centerline{Arnasli Yahya and Jen\H o Szirmai}}
{Optimal ball and horoball packings...}
\title
{Optimal ball and horoball packings generated by $3$-dimensional simply truncated Coxeter orthoschemes with parallel faces
\footnote{Mathematics Subject Classification 2010: 52C17; 52C22; 52B15. \newline
Key words and phrases: Coxeter group; horoball; hyperbolic geometry; packing; tiling \newline
}}

\author{Arnasli Yahya and Jen\H o Szirmai \\
\normalsize Budapest University of Technology and \\
\normalsize Economics Institute of Mathematics, \\
\normalsize Department of Geometry \\
\normalsize Budapest, P. O. Box: 91, H-1521 \\
\normalsize szirmai@math.bme.hu
\date{\normalsize{\today}}}
\maketitle
\begin{abstract}
In this paper we consider the ball and horoball packings belonging to $3$-dimensional Coxeter tilings that are 
derived by simply truncated orthoschemes with parallel faces.

The goal of this paper to determine the optimal ball and horoball packing arrangements and their densities
for all above Coxeter tilings in hyperbolic 3-space $\mathbb{H}^3$. 
The centers of horoballs are required to lie at ideal vertices of the
polyhedral cells constituting the tiling, and we allow horoballs of
different types at the various vertices. 

We prove that the densest packing of the above cases is realized by horoballs related to 
$\{\infty;3;6;\infty \}$ and $\{\infty;6;3;\infty \}$ tilings with density $\approx 0.8413392$.
\end{abstract}
\newtheorem{theorem}{Theorem}[section]
\newtheorem{corollary}[theorem]{Corollary}
\newtheorem{lemma}[theorem]{Lemma}
\newtheorem{exmple}[theorem]{Example}
\newtheorem{definition}[theorem]{Definition}
\newtheorem{rmrk}[theorem]{Remark}
\newtheorem{proposition}[theorem]{Proposition}
\newenvironment{remark}{\begin{rmrk}\normalfont}{\end{rmrk}}
\newenvironment{example}{\begin{exmple}\normalfont}{\end{exmple}}
\newenvironment{acknowledgement}{Acknowledgement}
%

%==================================================================%
%                             the main article                               %
%==================================================================%

%%%%%%%%%%%%%%%%%%%%%%%%%%%%%%%%%%%%%%%%%%%
%%%%%%%%%%%%%%%%%%%%%%%%%%%%%%%%%%%%%%%%%%%

%\bibliographystyle{ieeetr}
%
\section{Introduction, preliminary results}\label{intro}
In hyperbolic spaces $\mathbb{H}^n$ for $2 \leq n \leq 9$, the known densest ball and horoball configurations are derived by Coxeter simplex tilings, 
therefore it is interesting 
to determine their optimal horoball packings related to Coxeter tilings. In the former papers, we investigated that Coxeter simplex tilings 
whose generating simplices do not have parallel faces. Now, we extend our study to the $3$ dimensional Coxeter tilings generated by simple 
frustum orthoschemes with parallel faces. But first, we summarize the results related to the above topic.

Let $X$ denote a space of constant curvature, either the $n$-dimensional sphere $\mathbb{S}^n$, 
Euclidean space $\mathbb{E}^n$, or 
hyperbolic space $\mathbb{H}^n$ with $n \geq 2$. An important question of discrete geometry is to find the highest possible packing density in $X$ by 
congruent non-overlapping balls of a given radius \cite{Be}, \cite{G--K}. 
The definition of packing density is critical in hyperbolic space as shown by B\"or\"oczky \cite{BK1,BK2}, 
for the standard paradoxical construction see \cite{G--K}. 
The most widely accepted notion of packing density considers the local densities of balls with respect to their 
Dirichlet--Voronoi cells (cf. \cite{BK3} and \cite{KH1}). In order to study horoball packings in $\overline{\mathbb{H}}^n$, we use an extended notion of such local density. 

Let $B$ be a horoball of packing $\cB$, and $P \in \overline{\mathbb{H}}^n$ an arbitrary point. 
Define $d(P,B)$ to be the shortest distance from point $P$ to the horosphere 
$S = \partial B$, where $d(P,B)\leq 0$  if $P \in B$. 
The Dirichlet--Voronoi cell $\cD(B,\cB)$ of horoball $B$ is the convex body
\begin{equation}
\cD(B,\cB) = \{ P \in \mathbb{H}^n | d(P,B) \le d(P,B'), ~ \forall B' \in \cB \}. \notag
\end{equation}
Both $B$ and $\cD$ have infinite volume, so the standard notion of local density is
modified. Let $Q \in \partial{\mathbb{H}}^n$ denote the ideal center of $B$, and take its boundary $S$ to be the 
one-point compactification of Euclidean $(n-1)$-space.
Let $B_C^{n-1}(r) \subset S$ be the Euclidean $(n-1)$-ball with center $C \in S \setminus \{Q\}$.
Then $Q$ and $B_C^{n-1}(r)$ determine a convex cone 
$\cC^n(r) = cone_Q\left(B_C^{n-1}(r)\right) \in \overline{\mathbb{H}}^n$ with
apex $Q$ consisting of all hyperbolic geodesics passing through $B_C^{n-1}(r)$ with limit point $Q$. The local density $\delta_n(B, \cB)$ of $B$ to $\cD$ is defined as
\begin{equation}
\delta_n(\cB, B) =\varlimsup\limits_{r \rightarrow \infty} \frac{vol(B \cap \cC^n(r))} {vol(\cD \cap \cC^n(r))}. \notag
\end{equation}
This limit is independent of the choice of center $C$ for $B^{n-1}_C(r)$.

In the case of periodic ball or horoball packings, this local density defined above can be extended to the entire hyperbolic space, and 
is related to the simplicial density function (defined below) that we generalized in 
\cite{SzJ2} and \cite{SzJ3}.
In this paper, we shall use such definition of packing density (cf. Section 3).  

A Coxeter simplex is a top dimensional simplex in $X$ with dihedral angles either integral submultiples of $\pi$ or zero. 
The group generated by reflections on the sides of a Coxeter simplex is a Coxeter simplex reflection group. 
Such reflections generate a discrete group of isometries of $X$ with the Coxeter simplex as the fundamental domain; 
hence the groups give regular tessellations of $X$ if the fundamental simplex is characteristic. The Coxeter groups 
are finite for $\mathbb{S}^n$, and infinite for $\mathbb{E}^n$ or $\overline{\mathbb{H}}^n$. 

There are non-compact Coxeter simplices in $\overline{\mathbb{H}}^n$ with ideal vertices in $\partial \mathbb{H}^n$, 
however only for dimensions $2 \leq n \leq 9$; furthermore, only a finite number exists in dimensions $n \geq 3$. 
Johnson {\it et al.} \cite{JKRT} found the volumes of all mentioned Coxeter simplices in hyperbolic $n$-space, also see Kellerhals \cite{KH1}. 
Such simplices are the most elementary building blocks of hyperbolic manifolds,
the volume of which is an important topological invariant. 

In $n$-dimensional space $X$ of constant curvature
$(n \geq 2)$, define the simplicial density function $d_n(r)$ to be the density of $n+1$ mutually tangent balls of 
radius $r$ in the simplex spanned by their centers. L.~Fejes T\'oth and H.~S.~M.~Coxeter
conjectured that the packing density of balls of radius $r$ in $X$ cannot exceed $d_n(r)$.
Rogers \cite{Ro64} proved this conjecture in Euclidean space $\mathbb{E}^n$.
The $2$-dimensional spherical case was settled by L.~Fejes T\'oth \cite{FTL}, and 
B\"or\"oczky \cite{BK3}, who extend the analogous statement to 
$n$-dimensional sapces of constant curvature.

In hyperbolic 3-space, 
the monotonicity of $d_3(r)$ was proved by B\"or\"oczky and Florian
in \cite{BK4}; in \cite{Ma99} Marshall 
showed that for sufficiently large $n$, 
function $d_n(r)$ is strictly increasing in variable $r$. 

The simplicial packing density upper bound $d_3(\infty) = (1+\frac{1}{2^2}-\frac{1}{4^2}-\frac{1}{5^2}+\frac{1}{7^2}+\frac{1}{8^2}--++\dots)^{-1} = 0.85327\dots$ 
cannot be achieved by packing regular balls, instead it is realized by horoball packings of
$\overline{\mathbb{H}}^3$, the regular ideal simplex tiles $\overline{\mathbb{H}}^3$.
More precisely, the centers of horoballs in  $\partial\overline{\mathbb{H}}^3$ lie at the vertices of the ideal regular Coxeter simplex tiling with Schl\"afli symbol $\{3,3,6\}$. 

In \cite{KSz} we proved that this optimal horoball packing configuration in $\mathbb{H}^3$  is not unique. We gave several more examples of 
regular horoball packing arrangements based on asymptotic Coxeter tilings using horoballs of different types, that is horoballs that have different 
relative densities with respect to the fundamental domain, that yield the B\"or\"oczky--Florian-type simplicial upper bound \cite{BK4}.
 
Furthermore, in \cite{SzJ2,SzJ3} we found that 
by allowing horoballs of different types at each vertex of a totally asymptotic simplex and generalizing 
the simplicial density function to $\overline{\mathbb{H}}^n$ for $(n \ge 2)$,
the B\"or\"oczky-type density 
upper bound is not valid for the fully asymptotic simplices for $n \geq 4$. 
For example, in $\overline{\mathbb{H}}^4$ the locally optimal simplicial packing density is $0.77038\dots$, higher than the B\"or\"oczky-type density upper bound of $d_4(\infty) = 0.73046\dots$ using horoballs of a single type. 
However these ball packing configurations are only locally optimal and cannot be extended to the entirety of the
ambient space $\overline{\mathbb{H}}^n$. 
In \cite{KSz2} we found seven horoball packings of Coxeter simplex tilings in $\overline{\mathbb{H}}^4$ that 
yield densities of $5\sqrt{2}/\pi^2 \approx 0.71645$, counterexamples to L. Fejes T\'oth's conjecture of $\frac{5-\sqrt{5}}{4}$ stated in his foundational 
book {\it Regular Figures} \cite[p. 323]{FTL}.

In \cite{KSz3}, \cite{KSz4} we extend our study of horoball packings to $\overline{\mathbb{H}}^n$ ($5\le n \le 9$) 
using our methods that were successfully applied in lower dimensions.

In the previously mentioned papers, we studied the ball and horoball packing 
related to the Coxeter simplex tilings where the vertices of the simplices 
are proper points of the hyperbolic space 
$\mathbb{H}^n$ or they are ideal i.e. lying on 
the sphere $\partial \mathbb{H}^n$.

Now, we consider the Coxeter tilings in $3$-dimensional 
hyperbolic space $\mathbb{H}^3$ where the generating orthoscheme is a 
simply truncated Coxeter orthoscheme. 
Moreover, we require that the truncated orthoscheme has parallel 
faces i.e. their dihedral angle is zero. 

In $\mathbb{H}^3$ the simply truncated orthoscheme tilings (crystallographic 
Napier cycles of type 2 see \cite{IH1,IH2}) are determined by the following 
Coxeter-Schl\"afli graphs:
\begin{center}
\includegraphics[scale=0.3]{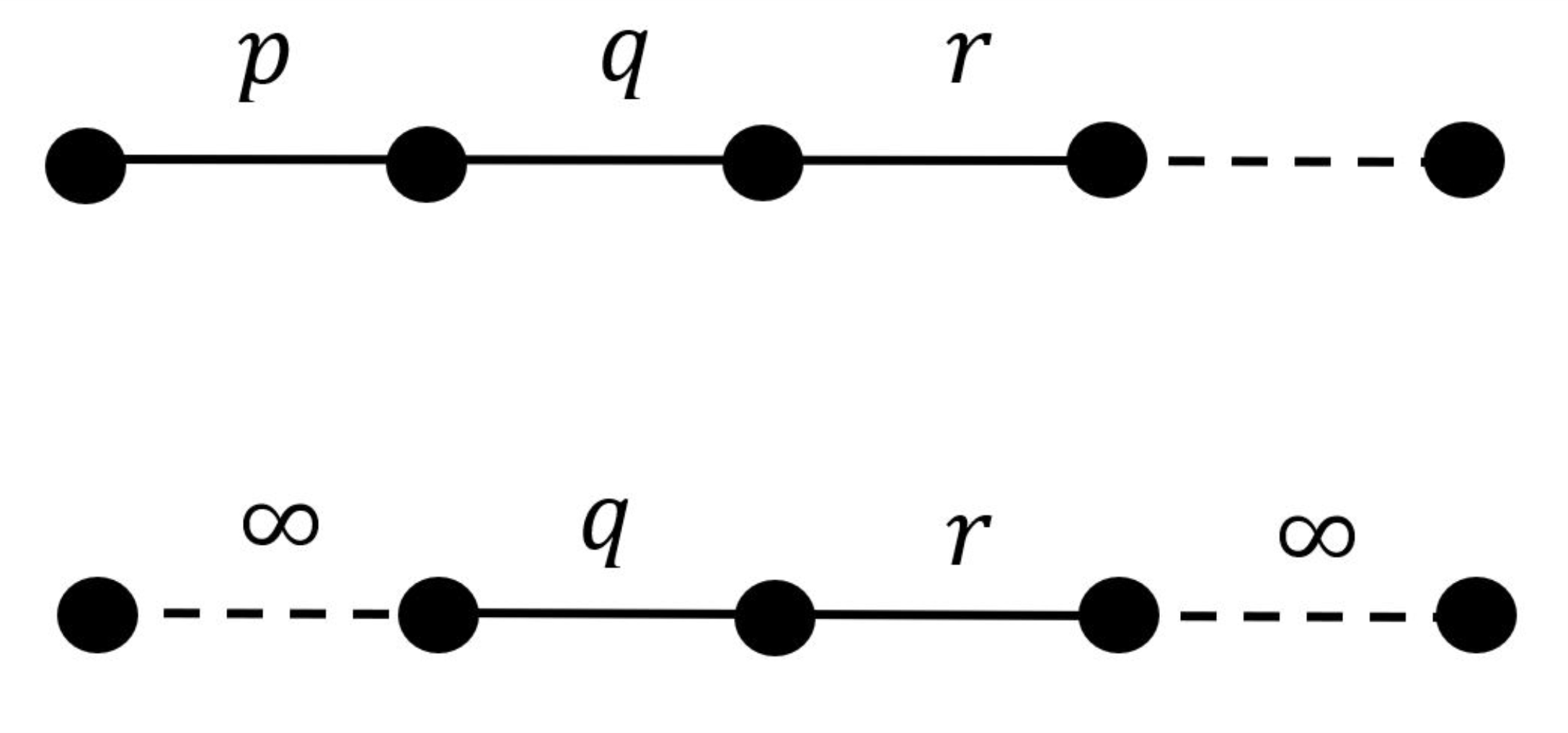}
\end{center}
where $\frac{1}{p}+\frac{1}{q} < \frac{1}{2}$ and $\frac{1}{q}+\frac{1}{r}\ge\frac{1}{2}$, ($3\leq p,q,r \in \mathbb{N}$).
\begin{rmrk}
The first cases are investigated in papers \cite{SzJ9,SzJ10,SzJ11} where from the 
truncated orthoschemes can be derived prism-like tilings generated hyperball 
or hyp-hor packings
(or coverings) in $\mathbb{H}^3$. 
\end{rmrk}
In this paper, we concentrate on the second tilings which are given with 
Schl\"afli symbol $\{\infty,q,r,\infty\}$. We determine the optimal ball and horoball 
packings
related to the above tilings prove that the densest packing arrangement of 
the considered tilings is realized at the tilings $\{\infty,3,6,\infty\}$, and $\{\infty;6;3;\infty \}$ 
by horoballs with density
$\approx 0.8413392$. Our results are summarized in Theorems 3.5, 4.2, 4.4 and in Tables 1,3,4.

\section{Basic notions}
\subsection{The projective model of hyperbolic space $\mathbb{H}^3$}
For the computations, we use the projective model \cite{ME} of the hyperbolic space. The model is defined in the $\mathbb{E}^{1,n}$ Lorentz space with signature $(1,n)$, i.e. consider $\mathbf{V}^{n+1}$ real vector space equipped with the bilinear form: 
\begin{equation*}
\langle ~ \mathbf{x},~\mathbf{y} \rangle = -x^0y^0+x^1y^1+ \dots + x^n y^n.
\end{equation*}
In the vector space, consider the following equivalence relation:
\begin{equation*}
\mathbf{x}(x^0,...,x^n)\thicksim \mathbf{y}(y^0,...,y^n)\Leftrightarrow \exists~ c\in \mathbb{R}\backslash\{0\}: \mathbf{y}=c\cdot\mathbf{x}
\end{equation*}
The factorization with $\thicksim$ induces the $\mathcal{P}^n(\mathbf{V}^{n+1},\mbox{\boldmath$V$}\!_{n+1})$ $n$-dimensional real projective space. In this space to 
interpret the points of $\mathbb{H}^n$ hyperbolic space, consider the following quadratic form:
\begin{equation*}
Q=\{[\mathbf{x}]\in\mathcal{P}^n | \langle ~ \mathbf{x},~\mathbf{x} \rangle =0 \}=:\partial \mathbb{H}^n
\end{equation*}
The inner points relative to the cone-component determined by $Q$ are the points of $\mathbb{H}^n$ (for them $\langle ~ \mathbf{x},~\mathbf{x} \rangle <0$), the point of $Q=\partial \mathbb{H}^n$ are called the points at infinity, and the points lying outside relative to $Q$ are outer points of $\mathbb{H}^n$ (for them $\langle ~ \mathbf{x},~\mathbf{x} \rangle >0$). We can also define a linear polarity between the points and hyperplanes of the space: the polar hyperplane of a point $[\mathbf{x}]\in\mathcal{P}^n$ is $Pol(\mathbf{x}):=\{[\mathbf{y}]\in\mathcal{P}^n | \langle ~ \mathbf{x},~\mathbf{y} \rangle =0 \}$, and hence $\mathbf{x}\in\mathbf{V}^{n+1}$ is incident with ${\boldsymbol{a}}\in\mbox{\boldmath$V$}\!_{n+1}$ iff $\langle ~ \mathbf{x},~{\boldsymbol{a}} \rangle =0$. In this projective model, we can define a metric structure related to the above bilinear form, where for the distance of two proper points:
\begin{equation}
\cosh(\frac{d(\mathbf{x},\mathbf{y})}{k})=\frac{-\langle ~ \mathbf{x},~\mathbf{y} \rangle}{\sqrt{\langle ~ \mathbf{x},~\mathbf{x} \rangle \langle ~ \mathbf{y},~\mathbf{y} \rangle}},
~ (\text{at present}~ k=1).
\end{equation}
This corresponds to the distance formula in the well-known Beltrami-Cayley-Klein model.
\subsection{Coxeter orthoschemes and tilings}
\begin{defn}
In the $\mathbb{H}^n$ $(2\leq n\in \mathbb{N})$ space a complete orthoscheme $\mathcal{O}$ of degree $d$ $(0\leq d\leq2)$ is a polytope bounded 
with hyperplanes $H^0,...,H^{n+d}$, for which $H^i\perp H^j$, unless $j\neq i-1,i,i+1$.
\end{defn}

In the classical ($d=0$) case, let denote the vertex opposite to $H^i$ hyperplane with $A_i$ ($0\leq i\leq n$), and let denote the dihedral angle of 
$H^i$ and $H^j$ planes with $\alpha^{ij}$ (hence $\alpha^{ij}$ if $0\leq i<j-1\leq n$). 

In this paper, we deal with orthoschemes of degree $d=1$, they can be described geometrically, as follows. 
We can give the sequence of the vertices of the orthoschemes $A_0,...,A_n$, where $A_iA_{i+1}$ edge is perpendicular 
to $A_{i+2}A_{i+3}$ edge for all $i\in\{0,...,n-3\}$. Here $A_0$ and $A_n$ are called the principal vertices of the orthoschemes. 
In the case $d=1$, one of these principal vertices (e.g. $A_0$) is the outer points of the above model, 
so they are truncated by its polar planes 
$Pol(A_0)$ and the orthoscheme is called simply truncated. 

In general, the Coxeter orthoschemes were classified by {\sc H-C. Im Hof}, he proved that they exist in dimension $\leq 9$, and gave a full list of them \cite{IH1}, \cite{IH2}. 

Now consider the reflections on the facets of the simly truncated orthoscheme, and denote them with  $r_1,...,r_{n+2}$, hence define the group 
\begin{equation*}
\mathbf{G}=\langle r_1,...,r_{n+2} | (r_ir_j)^{m_{ij}}=1 \rangle,
\end{equation*}
where $\alpha^{ij}=\frac{\pi}{m_{ij}}$, so $m_{ii}=1$, and if $m_{ij}=\infty$ (i.e. $H^i$ and $H^j$ are parallel), than to the $r_i,r_j$ pair belongs no relation. 
Suppose that $2\leq m_{ij}\in\{\mathbb{N}\cup\infty\}$ if $i\neq j$. The Coxeter group $\bG$ acts on hyperbolic 
space $\overline{\mathbb{H}}^n$ properly discontinously, thus the images of the orthoscheme under this action provide a 
$\mathcal{T}$ tiling of $\overline{\mathbb{H}}^n$ (i.e. the images of the orthoscheme fills $\overline{\mathbb{H}}^n$ without overlap).

For the {\it complete Coxeter orthoschemes} $\mathcal{S} \subset \mathbb{H}^n$, we adopt the usual
conventions and sometimes even use them in the Coxeter case: if two nodes are related by the weight $\cos{\frac{\pi}{m_{ij}}}$
then they are joined by a ($m_{ij}-2$)-fold line for $m_{ij}=3,~4$ and by a single line marked by $m_{ij}$ for $m_{ij} \geq 5$.
In the hyperbolic case if two bounding hyperplanes of $S$ are parallel, then the corresponding nodes
are joined by a line marked $\infty$. If they are divergent then their nodes are joined by a dotted line.

In the following we concentrate only on dimension $3$, $\overline{\mathbb{H}}^3$. 
For every considered tiling there is a corresponding
symmetric $4 \times 4$ matrix $(b^{ij})$ where $b^{ii}=1$ and, for $i \ne j\in \{0,1,2,3\}$,
$b^{ij}$ equals to $-\cos{\alpha_{ij}}$ with all dihedral angles $\alpha_{ij}$ 
between the faces $u_i$,$u_j$ of $\mathcal{S}$.

The matrix $(b^{ij})$ in formula (2.2) is the so called Coxeter-Schl\"afli matrix with
parameters $(p;q;r)$, i.e. $\alpha_{01}=\frac{\pi}{p}$, 
$\alpha_{12}=\frac{\pi}{q}$, $\alpha_{23}=\frac{\pi}{r}$.
Now only $3\le p,q,r \in \mathbb{N}$ come into account (see \cite{IH1,IH2}).
\begin{equation}
(b^{ij})=\langle \mbox{\boldmath$b^i$},\mbox{\boldmath$b^j$} \rangle:=\begin{pmatrix}
1& -\cos{(\frac{\pi}{p})}& 0 & 0 \\
-\cos{(\frac{\pi}{p})} & 1 & -\cos{(\frac{\pi}{q})}& 0 \\
0 & -\cos{(\frac{\pi}{q})} & 1 & -\cos{(\frac{\pi}{r})} \\
0 & 0 & -\cos{(\frac{\pi}{r})} & 1
\end{pmatrix}. 
\end{equation}

This $3$-dimensional complete (truncated or frustum) orthoscheme $\mathcal{S}=\mathcal{S}(p;q;r)$ tilings are described in
Fig.~5,~6,~7 and they are characterised by their symmetric Coxeter-Schl\"afli matrices $(b^{ij})$ (see formula (2.2)), furthermore by their inverse matrices $(h_{ij})$ 
in formula (2.3).
\begin{equation}
\begin{gathered}
(h_{ij})=(b^{ij})^{-1}=\langle \ba_i, \ba_j \rangle:=\\
=\frac{1}{B} \begin{pmatrix}
\sin^2{(\frac{\pi}{r})}-\cos^2{(\frac{\pi}{q})}& \cos{(\frac{\pi}{p})}\sin^2{(\frac{\pi}{r})}& \cos{(\frac{\pi}{p})}\cos{(\frac{\pi}{q})} & \cos{(\frac{\pi}{p})}\cos{(\frac{\pi}{q})}\cos{(\frac{\pi}{r})} \\
\cos{(\frac{\pi}{p})}\sin^2{(\frac{\pi}{r})} & \sin^2{(\frac{\pi}{r})} & \cos{(\frac{\pi}{q})}& \cos{(\frac{\pi}{r})}\cos{(\frac{\pi}{q})} \\
\cos{(\frac{\pi}{p})}\cos{(\frac{\pi}{q})} & \cos{(\frac{\pi}{q})} & \sin^2{(\frac{\pi}{p})}  & \cos{(\frac{\pi}{r})}\sin^2{(\frac{\pi}{p})}  \\
\cos{(\frac{\pi}{p})}\cos{(\frac{\pi}{q})}\cos{(\frac{\pi}{r})}  & \cos{(\frac{\pi}{r})}\cos{(\frac{\pi}{q})} & \cos{(\frac{\pi}{r})}\sin^2{(\frac{\pi}{p})}  & \sin^2{(\frac{\pi}{p})}-\cos^2{(\frac{\pi}{q})}
\end{pmatrix}, 
\end{gathered}
\end{equation}
where
$$
B=\det(b^{ij})=\sin^2{\left(\frac{\pi}{p}\right)}\sin^2{\left(\frac{\pi}{r}\right)}-\cos^2{\left(\frac{\pi}{q}\right)} <0, \ \ \text{i.e.} \ \sin{\left(\frac{\pi}{p}\right)}\sin{\left(\frac{\pi}{r}\right)}-\cos{\left(\frac{\pi}{q}\right)}<0.
$$
In this work, our aim is to study the Coxeter tilings generated by simply truncated orthoschemes with Coxeter-Schl\"afli graph $\{\infty;q;r;\infty\}$. In these cases, 
the truncated orthoscheme is a $3$-dimensional hyperbolic polyhedron bounded by $5$ faces $u_i$ ($i=0,\dots 4$) that are determined by its form $\Bu_i$.
If the parameters $p, q, r$ satisfy $\frac{1}{p}+\frac{1}{q}<\frac{1}{2}$ and $\frac{1}{q}+\frac{1}{r} \geq \frac{1}{2}$,  ($3\le p,q,r \in \mathbb{N})$ then the 
corresponding Coxeter-Schl\"afli matrix is
\begin{equation} \label{CoxeterMatrix}
    (B^{ij})=\begin{pmatrix}
    1 & -\cos{(\frac{\pi}{p})}&0&0&0\\
    -\cos{(\frac{\pi}{p})}&1&-\cos{(\frac{\pi}{q})}&0&0\\
    0&-\cos{(\frac{\pi}{q})}&1&-\cos{(\frac{\pi}{r})}&0\\
    0&0&-\cos{(\frac{\pi}{r})}&1&c_4\\
    0&0&0&c_4&1
    \end{pmatrix},
\end{equation}
where the constant $c_4$ can be uniquely determined in the arrangement of Napier cycles 
\cite{IH1}\\
\begin{equation*}
    c_4 = -\sqrt{\frac{1+\cos^2(\frac{\pi}{p})\cos^2(\frac{\pi}{r})-\cos^2(\frac{\pi}{p})-\cos^2(\frac{\pi}{q})-\cos^2(\frac{\pi}{r})}{1-\cos^2(\frac{\pi}{p})-\cos^2(\frac{\pi}{q})}}.
\end{equation*}
In our case, there are two parallel faces (hyperplanes) that meet at the ideal point. 
That means the dihedral angle between these two hyperplanes are equal to zero. 
Therefore, we assume that these two hyperplane are $u_0$ and $u_1$, in such their dihedral angle is $w_0=\frac{\pi}{p}=0$. 
The corresponding Coxeter-Schl\"afli matrix $(B^{ij})$ can be described in the following form:
\begin{equation}\label{CexeterInfinity}
    (B^{ij})=\begin{pmatrix} 
    1 & -1&0&0&0\\
    -1&1&-\cos{(\frac{\pi}{q})}&0&0\\
    0&-\cos{(\frac{\pi}{q})}&1&-\cos{(\frac{\pi}{r})}&0\\
    0&0&-\cos{(\frac{\pi}{r})}&1&-1\\
    0&0&0&-1&1
    \end{pmatrix}
\end{equation}
The matrix $(B^{ij})$ is a singular matrix whose non-singular $4 \times 4$ principal submatrix is $(b^{ij})$.

The volume of a simply truncated Coxeter orthoscheme with outer vertices $A_0$ can be determined by the 
following theorem of {\sc R. Kellerhals} \cite{KH1,KH2}.
\begin{theorem} The volume of a three-dimensional hyperbolic
complete ortho\-scheme (except Lambert cube cases) $\mathcal{S}$
is expressed with the essential angles $\alpha^{01},\alpha^{12},\alpha^{23}, \\ (0 \le \alpha^{ij} \le \frac{\pi}{2})$ in the following form:

\begin{align*}
&Vol_3(\mathcal{S})=\frac{1}{4} \{ \mathcal{L}(\alpha^{01}+\theta)-
\mathcal{L}(\alpha^{01}-\theta)+\mathcal{L}(\frac{\pi}{2}+\alpha^{12}-\theta)+ \notag \\
&+\mathcal{L}(\frac{\pi}{2}-\alpha^{12}-\theta)+\mathcal{L}(\alpha^{23}+\theta)-
\mathcal{L}(\alpha^{23}-\theta)+2\mathcal{L}(\frac{\pi}{2}-\theta) \}, 
\end{align*}
where $\theta \in [0,\frac{\pi}{2})$ is defined by the following formula:
$$
\tan(\theta)=\frac{\sqrt{ \cos^2{\alpha^{12}}-\sin^2{\alpha^{01}} \sin^2{\alpha^{23}
}}} {\cos{\alpha^{01}}\cos{\alpha^{23}}}
$$
and where $\mathcal{L}(x):=-\int\limits_0^x \log \vert {2\sin{t}} \vert dt$ \ denotes the
Lobachevsky function.
\end{theorem}

In this paper, we consider that Coxeter orthoscheme tilings which are given with Schl\"afli symbol $\{\infty;q;r;\infty\}$ ($\frac{1}{q}+\frac{1}{r} \ge \frac{1}{2}$, 
$3 \leq p,q,r \in \mathbb{N}$) and determine the optimal ball and horoball packings related to the above tilings. 
The possible parameters are the following:
\begin{equation}
(q,r)=(3,4),(3,5),(3,6),(4,3),(4,4),(4,5),(5,3),(6,3).
\end{equation}
\section{Optimal ball packings}
\subsection{Inradii of truncated orthoschemes}
In determining the radius entered, we must distinguish two cases. The first type includes cases where the inscribed 
ball of the complete orthoscheme is the same as the inscribed ball of the truncated orthoscheme, and the second case where it is not true. 
\subsubsection{Type 1}
In \cite{J14} {\sc M. Jacquemet} determined the inradii of truncated simplices in $n$-dimensi-onal hyperbolic spaces if their inballs do not have common inner points with the 
corresponding truncating hyperplanes. We first recall some of the statements from that mentioned paper that we will apply to our calculations. 

We use the former section introduced denotation: let $\mathcal{S}$ a complete orthoscheme and $\hat{\mathcal{S}}$ be  the considered truncated orthoscheme. 
The corresponding Coxeter-Schl\"afli matrix is denoted by $(B^{ij})$ whose principal submatrix $(b^{ij})$.

\begin{lemma}\label{inradiusConditions}
A truncated hyperbolic simplex $\hat{\mathcal{S}}$ with Coxeter-Schl\"afli principal submatrix $(b^{ij})$ has inball, (imbedded ball of maximal finite radius) in $\mathbb{H}^n$ if only 
if $\sum_{i,j=1}^{n+1} cof_{ij}({b^{ij}}) > 0$.
\end{lemma}
We can rewrite this conditions using $h_{ij}$ (the inverse matrix of $(b^{ij}))$.
\begin{equation}
    \sum_{i,j=1}^{n+1} det(b^{ij})(h_{ij}) >0
\end{equation}
\begin{lemma}\label{Inradius}
Let $(b^{ij})$ be the Coxeter-Schl\"afli matrix of complete orthoscheme ${\mathcal{S}}$ with inball ${\mathcal{B}} \subset \mathbb{H}^n$. 
Then, the inradius ${r}=r({\mathcal{B}})$ is given by
\begin{equation}
{r}=\sinh^{-1}{\sqrt{-\frac{1}{\sum_{i,j=1}^{n+1} h_{ij}}}}.
\end{equation}
\end{lemma}
In our case, the orthoscheme is truncated by a polar hyperplane of the vertex which lie outer of the Beltrami-Cayley-Klein model, 
it is called ultra ideal vertex. Therefore, we need to study whether the inradius of the above orthoscheme and the inradius of the truncated orthoscheme are equal to each other.
\begin{lemma}[\cite{J14}]
Let ${\mathcal{S}}$ be an $n$-dimensional simplex with $k \le n$ vertices are ultra ideal, with Coxeter-Schl\"afli matrix $(b^{ij})$, such that ${\mathcal{S}}$ has an inball 
${B} \subset \mathbb{H}^n$ of radius ${r}$. Denote by $\hat{\mathcal{S}} \subset \bar{\mathcal{H}^n}$ its associated hyperbolic $k$-truncated simplex 
with respect to the ultra-ideal vertices $v_1, \cdots, v_k$, $1 \leq k$. Let $\hat{r}$ be the inradius of inball of $\hat{\mathcal{S}}$. Then, $r=\hat{r}$ if and only if
\begin{equation}
    \frac{\sum_{j=1}^{n+1}cof_{ij}{(b^{ij})}}{det({b^{ij}})cof_{ii}{(b^{ij})}} \geq 1 ~~~ \text{for all}~ i=1, \cdots, k.
\end{equation}
\end{lemma}
\begin{figure}[ht]
\centering
\includegraphics[width=\textwidth]{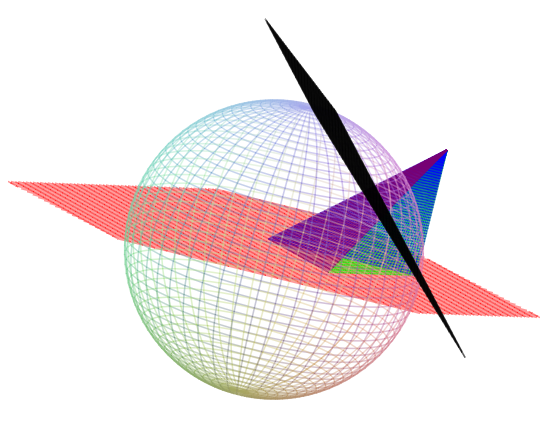}
\caption{The truncated orthoscheme: $(q,r)=(3,4)$}
\end{figure}
\begin{figure}[ht]
\centering
\includegraphics[width=8 cm]{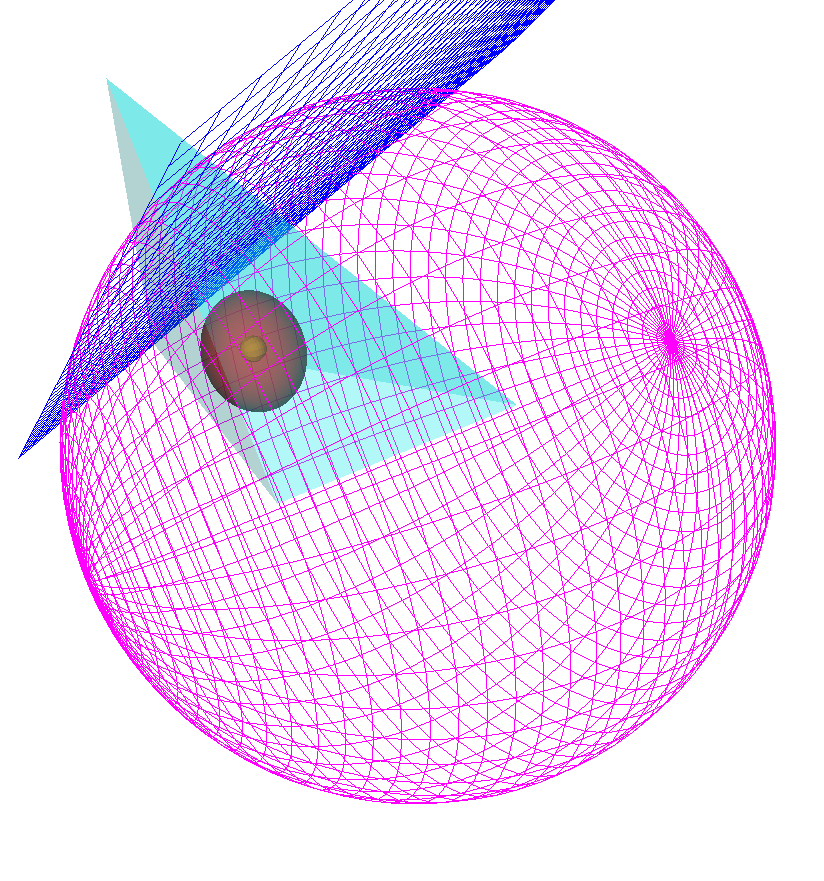}
\caption{The structure of truncated orthoscheme with inball and its centre}
\end{figure}

We apply the above existence conditions for the considered tilings and obtain, ${\det}{(b^{ij})}=- \cos^{2}{\left(\frac{\pi}{q} \right)}$.
To confirm the existence of inball in the complete orthoscheme $\mathcal{S}$ we need to clarify whether the value 
$\sum_{i,j=1}^{n+1} det(b^{ij})h_{ij}$ is a positive number. By the direct computations, we have the following expression
{\footnotesize{
\begin{equation}
\sum_{i,j=1}^{n+1} det(b^{ij})h_{ij}=  \left( \frac{4+2\cos^{2}\left( \frac{\pi}{q}\right)+
4\cos\left( \frac{\pi}{q}\right)\cos\left( \frac{\pi}{r}\right)-
4 \cos^{2}\left( \frac{\pi}{r}\right)}{\cos^{2}\left( \frac{\pi}{q}\right)} \right) \cos^{2} \left( \frac{\pi}{q} \right). 
\end{equation}}}
One could see that in the expression above, 
\begin{equation}
\begin{gathered}
4+2\cos^{2}\left( \frac{\pi}{q}\right)+4\cos\left( \frac{\pi}{q}\right)\cos\left( \frac{\pi}{r}\right)-4 \cos^{2}\left( \frac{\pi}{r}\right)> \\
> 4-4\cos^{2}\left( \frac{\pi}{r}\right)> 0.
\end{gathered}
\end{equation}
Therefore $\sum_{i,j=1}^{n+1} det(b^{ij})h_{ij}>0 $ for all $(q,r)$ in our cases. 

Then, we have to determine for which parameters are the inradius $r$ of complete orthoscheme $\mathcal{S}$ equal to the inradius $\hat{r}$ of the truncated orthoscheme 
$\hat{\mathcal{S}}$.

In our case there is only one ultra ideal vertex therefore, applying the formula (3.3) we get the following inequalities:
\begin{equation}
\begin{gathered}
 - \frac{\left(1 - \frac{2 \cos{\left(\frac{\pi}{r} \right)}}{\cos{\left(\frac{\pi}{q} \right)}}\right) 
 \cos^{2}{\left(\frac{\pi}{q} \right)}}{\sqrt{\cos^{4}{\left(\frac{\pi}{q} \right)}}} \geq 1 \Leftrightarrow 
  \frac{\cos{\left(\frac{\pi}{r} \right)}}{\cos{\left(\frac{\pi}{q} \right)}} \geq 1.
\end{gathered}
\end{equation}

That inequality could only be satisfied by the parameter tuples of $(q,r)=$ $(3,3)$, $(3,4)$, $(3,5)$, $(3,6)$, $(4,4)$.\\
Therefore, in these cases the radii of optimal inballs can be determined by formula (3.2).
\subsection{Type 2}
We consider the case where the inradius of the complete orthoschemes and the truncated orthoschemes are not the same, 
i.e when the Lemma 3.3 does not hold. This situation happens on the tuple parameters 
$(4,3),(5,3),(6,3)$. In these cases, the constructed inball intersect the truncated face:
\begin{align*}
r=d(\Bc,\Bu_0)=d(\Bc,\Bu_1)=d(\Bc,\Bu_2)=d(\Bc,\Bu_3) > d(\Bc,\Bu_4),
\end{align*}
where $\Bc$ is the insphere center of the complete orthoscheme $\mathcal{S}(q,r)$. 
Therefore, the problem is that the constructed inball intersects the hyperplane $\Bu_4$.

Instead of using the results of \cite{J14}, we apply the classical way to determine the incenter and its radius in the mentioned cases.
Basically, to find the center of optimal inball we determine the hyperplane bisectors of faces of truncated orthoscheme $\hat{\mathcal{S}}(q,r)$. 
The inball of maximal radius has to touch at least four faces of $\hat{\mathcal{S}}(q,r)$ and one of them must be the face determined by form $[\Bu_4]$. 

Therefore, we have $4$ analogues cases that provide the candidates for the optimal incenter. We have to determine these centers and 
we need to select the center with the maximum radius. This optimal ball is denoted by $\mathcal{B}(q,r)$. 

\subsection{Densest ball packing configurations}
The volume of the truncated orthoscheme $\hat{\mathcal{S}}(q,r)$ is denoted by $Vol(\hat{\mathcal{S}}(q,r))$ 
where $(q,r)=$ $(3,3)$, $(3,4)$, $(3,5)$, $(3,6)$, $(4,4)$, $(4,3)$, $(5,3)$, $(6,3)$.
We introduce the local density function $\delta_{opt}(\hat{\mathcal{S}(q,r)})$ related to orthoscheme $\hat{\mathcal{S}(q,r)}$ generated tiling:
\begin{definition}
The local density function $\delta_{opt}(\hat{\mathcal{S}}(q,r))$ related to tiling $\cT(q,r)$ generaded by truncated orthoscheme $\hat{\mathcal{S}}(q,r)$:
\begin{equation}
\delta_{opt}(\hat{\mathcal{S}}(q,r)):=\frac{Vol(\mathcal{B}(q,r))}{Vol(\hat{\mathcal{S}}(q,r))}. \notag
\end{equation}
\end{definition}
We obtain the volumes of the $Vol(\mathcal{B}(q,r))$ balls by the classical formula of {\sc J. Bolyai} 
$Vol(\mathcal{B}(q,r))=\pi k^3\left( \sinh{\left( \frac{2r}{k} \right)}-\frac{2r}{k}\right)$, where at present $k=1$ and 
the volumes of truncated simplices 
$\hat{\mathcal{S}}(q,r)$ can be calculated by Theorem 2.1. Using the formulas of the previous subsection we summarize our results in the following
table:
\begin{center}
Table 1.

 \vspace{2mm}
 \begin{tabular}{||c c c c c||} 
 \hline
 $(q,r)$ & Inradius $\hat{r}$ & $Vol(\mathcal{B}(q,r))$ & $Vol(\hat{\mathcal{S}}(p,q))$ & $\delta_{opt}(\hat{\mathcal{S}}(q,r))$ \\ [0.5ex] 
 \hline\hline
 $(3,3)$ & 0.2116177  & 0.0400529 & 0.1526609 & 0.2623649 \\ 
 \hline
 $(3,4)$ & 0.2236802 & 0.0473496 & 0.2509603 & 0.1886735\\
 \hline
 $(3,5)$ & 0.2335727 & 0.0539625 & 0.3323272 & 0.1623776\\
 \hline
 $(3,6)$ & 0.2407179 & 0.0591079 & 0.4228923 & 0.1397706 \\
 \hline
 $(4,3)$ &0.2396177&0.0582950&0.2509603&0.2322876 \\
 \hline
 $(4,4)$ & 0.2888593 & 0.1026579 & 0.4579828 & 0.2241524 \\ 
 \hline
 $(5,3)$ &0.2562904&0.0714478&0.3323273&0.2149924 \\
 \hline
 $(6,3)$ &0.2431555&0.0609361&0.4228923&0.1440937\\
 \hline
\end{tabular}
\end{center}
Finally, we obtain the following
\begin{theorem}
In hyperbolic space $\mathbb{H}^3$, between congruent ball packings of {\it classical balls}, generated by simply truncated Coxeter orthoschemes with parallel faces, 
the $\mathcal{B}{(3,3)}$ ball configuration provides the densest packing with density 
$\approx 0.2623649$. 
\end{theorem}
\section{Horoball packings}
The aim of this section is to determine the optimal horoball packing
densities for Coxeter tilings which are given by Schl\"afli symbol $\{\infty,q,r,\infty\}$ in 3-dimensional 
hyperbolic space $\mathbb{H}^3$ and have at least one vertex at the infinity $\partial \mathbb{H}^3$
$(q,r)=$ $(3,3)$, $(3,4)$, $(3,5)$, $(3,6)$, $(4,4)$, $(4,3)$, $(5,3)$, $(6,3)$ their truncated orthoschemes 
(fundamental domains of of the corresponding Coxeter group $\Gamma(q,r)$). The volume of the truncated orthoscheme 
$\hat{\mathcal{S}}(q,r)$ is denoted, similarly to the former subsection, by $Vol(\hat{\mathcal{S}}(q,r))$.
Moreover, $\mathcal{T}(q,r)$ denotes the corresponding hyperbolic truncated orthoscheme tiling.   
In cases $(q,r)=(3,6),(4,4), (6,3)$ the fundamental polyhedron has two vertices 
lying on the absolute $\partial \mathbb{H}^3$.  

An approach to describing Coxeter tilings involves analysis of their symmetry
groups. If $\cT(q,r)$ is a Coxeter orthoscheme tiling, then any rigid motion moving one cell into
another maps the entire tiling onto itself. 
Any orthoscheme cell of $\cT(q,r)$ can act as the fundamental domain $\hat{\mathcal{S}}(q,r)$
of $\Gamma(q,r)$ generated by reflections on its $2$-dimensional faces. 

We define the density of a horoball packing $\mathcal{B}(q,r)$ related to Coxeter orthoscheme tiling $\cT(q,r)$ as
\begin{equation}
\delta(\mathcal{B}({q,r}))=\frac{\sum_{i=1}^k Vol(\mathcal{B}_i \cap \hat{\mathcal{S}}(q,r))}{Vol(\hat{\mathcal{S}}(q,r))}
\end{equation}
where $\hat{\mathcal{S}}(q,r)$ denotes the fundamental domain orthoscheme of the tiling $\cT(q,r)$, $k\in\{1,2\}$ is the number of ideal vertices of $\hat{\mathcal{S}}(q,r)$, 
and $\mathcal{B}_i$ are the horoballs centered at the ideal vertices. 
We allow horoballs of different types at the ideal vertices of the tiling. 
A horoball type is allowed if it gives a packing, i.e. no two horoballs have an interior point in common. 
In addition, we require that no horoball may extend beyond the facet opposite the vertex where it is centered in order that the packing preserves the Coxeter symmetry group of the tiling.
If these conditions are satisfied we can use the corresponding Coxeter group $\Gamma(q,r)$ associated to a tiling to extend the packing density from 
the fundamental domain simplex to all of $\mathbb{H}^3$. We denote the optimal horoball packing density 
\begin{equation}
\delta_{opt}(\cB(q,r)) = \sup\limits_{\mathcal{B}(q,r) \text{~packing}} \delta(\mathcal{B}(q,r)).
\end{equation}
\subsubsection{Equations of the horospheres}
A horoshere in the hyperbolic geometry is the surface orthogonal to
the set of parallel lines, passing through the same point on the absolute quadratic surface $\partial \mathbb{H}^n$ (at present $n=3$)
(simply absolute) of  the hyperbolic space.

We represent hyperbolic space $\mathbb{H}^3$ in the Beltrami-Cayley-Klein
ball model. We introduce a projective coordinate system using vector
basis $\boldsymbol{b}_i \ (i=0,1,2,3)$ for $\mathcal{P}^3$ where the
coordinates of center of the model is $C=(1,0,0,0)$. We pick an
arbitrary point at infinity to be $A_3=(1,0,0,1)$.

As it is known, the equation of a horosphere with center
$A_3=(1,0,0,1)$ through point $S=(1,0,0,s)$ is 
\begin{gather}
\text{If $s \neq 1$, then} \ \ \  \frac{{(s-1)}^2}{1-s^2}(-x^0 x^0 +x^1 x^1+x^2 x^2+x^3 x^3)+{(x^0-x^3)}^2=0 \Leftrightarrow \notag \\
\Leftrightarrow (s-1)(-x^0 x^0 +x^1 x^1+x^2 x^2+x^3 x^3)-(1+s){(x^0-x^3)}^2=0 \notag
\end{gather}
Therefore, we obtain the following equation for the horosphere in 
Beltrami-Cayley-Klein model related to our Cartesian
coordinate system:
($x:=\frac{x^1}{x^0},~y:=\frac{x^2}{x^0},~z:=\frac{x^3}{x^0}$)
\begin{equation}
\frac{2(x^2+y^2)}{1-s}+\frac{4(z-\frac{s+1}{2})^2}{{(1-s)}^2}=1  
\end{equation}
\subsubsection{Volumes of horoball sectors}
The length $L(x)$ of a horospheric arc of a chord segment $x$ is determined
by the classical formula due to {\sc{J.~Bolyai}}:
\begin{equation}
L(x)=2 \cdot k \sinh{\frac{x}{2k}} \ (\text{at present} \ k=1).  
\end{equation}

The intrinsic geometry of the horosphere is Euclidean, therefore,
the area $\mathrm{Area}(A)$ of a horospherical triangle $A$ is computed by
the formula of Heron or by Cayley-Menger determinant. 
The volume of the horoball pieces can be
calculated using another formula by {\sc{J.~Bolyai}}. If the area of a domain on the horoshere is $\mathrm{Area}(A)$, then 
the volume determined
by $A$ and the aggregate of axes drawn from $A$ is equal to
\begin{equation}
V=\frac{1}{2}k \cdot \mathrm{Area}(A) \ \ (\text{at present} \ k=1). 
\end{equation}
\subsection{Horoball packings with one horoball}

We compute optimal horoball packing density for the Coxeter simplex tiling in case $(q,r)=(3,3)$. 
The other cases are obtained by using the same method. 

\begin{proposition}
\label{proposition:s4}
The optimal horoball packing density for Coxeter orthoscheme tiling of 
Schl\"afli symbol $\{\infty,3,3,\infty\}$
is $\delta_{opt}(\mathcal{B}(3,3)) \approx 0.8188080$.
\end{proposition}
{\it Proof:}
For the fundamental polyhedron (truncated orthoscheme) $\hat{\mathcal{S}}(3,3))$ of the Coxeter tiling $\cT(3,3)$ 
we fix coordinates for the vertices $A_0, A_1, A_2, A_4, A_5$ (see Fig.~3) 
to satisfy the angle requirements. 
The faces of $\hat{\mathcal{S}}(3,3))$ are given by their forms 
$[\mbox{\boldmath$u$}_i]$ $(i=0,1,2,3,4\})$ where $[\mbox{\boldmath$u$}_4]$ is the polar plane of outer point $A_3$ (see Fig.~3) and 
the other planes corresponding to the opposite vertex $A_i$ $(i=0,1,2,3)$ related to 
complete orthoscheme ${\mathcal{S}}(3,3))$.
\begin{figure}[ht!]
 \begin{center}
    \includegraphics[scale=1]{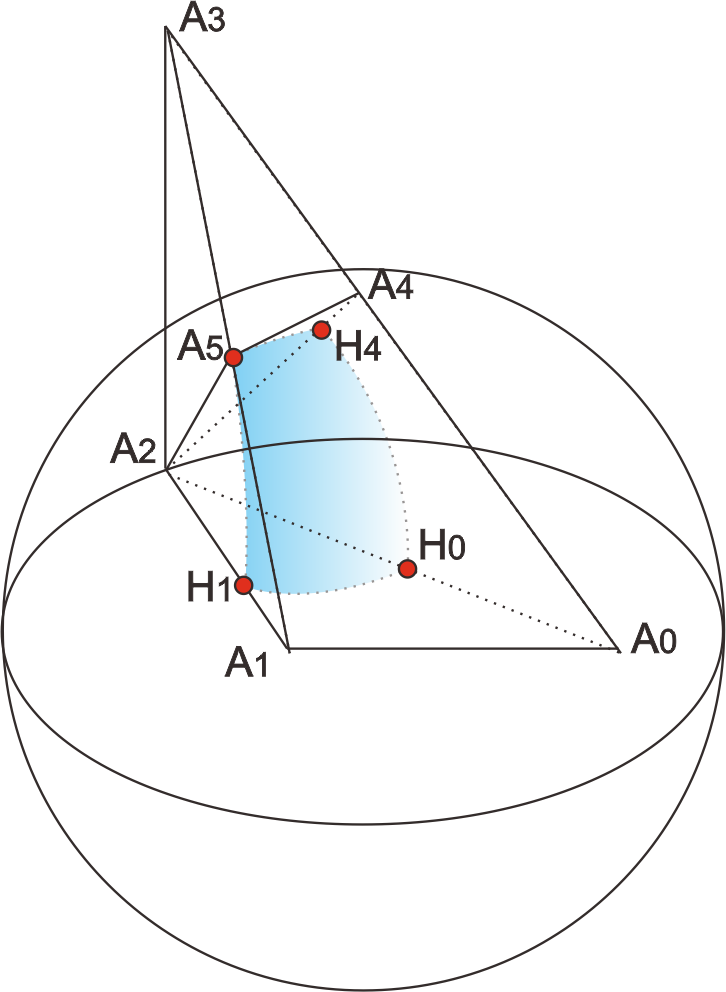}
    \vspace{0.1 cm}
    \caption{Optimal packing with one horoball type}
    \label{Footpoints}
    \end{center}
\end{figure}

To maximize the packing density, we determine the maximal horoball 
type $\mathcal{B}_2(s)$ centered at ideal vertex $A_2$ that fits into 
the fundamental domain $\hat{\mathcal{S}}(3,3))$. We find the horoball type parameter 
$s$ corresponding to the {\it radius} of the horoball when the horoball $\mathcal{B}_2(s)$ is tangent to 
the hyperface plane $[\mbox{\boldmath$u$}_2]$ bounding 
the fundamental polyhedron opposite of $A_2$. 
The perpendicular foot $F_2[\mbox{\boldmath$f$}_2]=A_5[\mathbf{a}_5]$ of vertex $A_2$ on 
plane $[\mbox{\boldmath$u$}_2]$,
\begin{equation}
\mathbf{a}_5 =\ba_2 - \frac{\langle \ba_2, \mathbf{u}_2 \rangle}{\langle \mathbf{u}_2,\mathbf{u}_2 \rangle} 
\mathbf{u}_2 \approx \left(1,~-0.074,~-0.430, 0.257 \right),
\end{equation}
is the point of tangency of the maximal horoball and face $[\mbox{\boldmath$u$}_2]$ of the orthoscheme cell.

Plugging in $F_2=H_5$ and solving equation (4.3) we get the parameter of the optimal horoball type. 
The equation of horosphere $\mathcal{B}_2$ centered at $A_2$ passing through $F_2=A_5$ can be determined by formulas (4.3) and (4.5).

The intersections $H_i[\bh_i]$ of the horosphere $\mathcal{B}_2$ and the orthoscheme edges are found by parameterizing the 
simplex edges as $\bh_i(\lambda) = \lambda \ba_2+\ba_i$ $(i=1,2,3,4)$, and computing their intersections with the horopshere (see Fig.~3).

The volume of the horospherical quadrilateral $A=H_0H_1H_4H_5$ determines the volume of the horoball piece by equation (4.5).
In order to determine the data of the horospheric quadrilateral we compute the hyperbolic distances $l_{ij}$ by the formula (2.1)
$l_{ij} = d(H_i, H_j)$. Moreover, the horospherical distances $L_{ij}$ can be calculated by the formula (4.4).

In the following table, we summarize the corresponding hyperbolic distances, horospherical distances: 
\begin{center}
Table 2.

\vspace{2mm}
 \begin{tabular}{||c c c||}
 \hline
 Edge $H_iH_j$  & Hyperbolic distance $l_{ij}$ & Horospherical distance $L_{ij}$  \\ [0.5ex] 
 \hline\hline
 $H_0 H_1$ & 0.4949329 & 0.5000000  \\ 
 \hline
 $H_1 H_4$ & 0.4949329 & 0.5000000 \\
 \hline
 $H_0 H_4$ & 0.6931471 & 0.7071067 \\
 \hline
 $H_4 H_5$ & 0.4949329 & 0.5000000 \\
 \hline
 $H_0 H_5$ & 0.4949329 & 0.5000000 \\[1ex] 
 \hline
\end{tabular}
\end{center}
To determine the area of quadrilateral $A=H_0H_1H_4H_5$, divide it into two horospheric triangles using a horospheric diagonal curve.
The intrinsic geometry of the horosphere is Euclidean so we can use it to find the area $\mathrm{Area}(A)$ of $A$ the 
Cayley-Menger determinant:
\begin{equation}
\begin{gathered}
    \mathrm{Area}(H_0H_1H_4)=\left(\frac{(-1)^{3}}{(2!)^{2} 2^{2}} 
    \begin{vmatrix} 0 & 1 & 1 & 1 \\ 1 & 0 & L_{0,1}^{2} & L_{04}^2 \\ 1 & L_{10}^2 & 0 & L_{14}^2 \\ 1 & L_{40}^2 & L_{41}^2 & 0 \end{vmatrix}\right)^{\frac{1}{2}}
    \approx0.1250000,\\ 
    \mathrm{Area}(H_0H_4H_5)=\left(\frac{(-1)^{3}}{(2!)^{2} 2^{2}} 
    \begin{vmatrix} 0 & 1 & 1 & 1 \\ 1 & 0 & L_{04}^{2} & L_{05}^2 \\ 1 & L_{40}^2 & 0 & L_{45}^2 \\ 1 & L_{50}^2 & L_{54}^2 & 0 \end{vmatrix}\right)^{\frac{1}{2}}
    \approx0.1250000, \\
    \Rightarrow \mathrm{Area}(A)=\mathrm{Area}(H_0H_1H_4)+\mathrm{Area}(H_0H_4H_5) \approx 0.2500000.
\end{gathered}
\end{equation}
The volume of the optimal horoball piece contained in the fundamental truncated orthoscheme $\hat{\mathcal{S}}(3,3))$ is

\begin{equation}
Vol(\mathcal{B}_2 \cap \hat{\mathcal{S}}(3,3))) = \frac{1}{2}\mathrm{Area}(A) \approx 0.1250000.
\end{equation}
Hence, the optimal horoball packing density of the Coxeter orthoscheme tiling $\mathcal{T}(3,3)$ becomes
\begin{equation}
\delta_{opt}(\mathcal{B}({3,3}))=\frac{Vol(\mathcal{B}_2 \cap \hat{\mathcal{S}}(3,3))}{Vol(\hat{\mathcal{S}}(q,r))}\approx 0.8188080. 
\end{equation}
The same method is used to find the optimal packing density of the remaining Coxeter truncated orthoscheme tilings. 
In cases if the truncated orthoscheme has two ideal vertices (vertices lying at the infinity)
then the optimal horoball packing configurations with one horoball are interesting, too. Therefore, we determine both optimal packing arrangements and their densities.

The results of the computations are summarized in Table 3 (see Fig.~4, 5).
\begin{center}
Table 3. 

\vspace{2mm}
 \begin{tabular}{||c c c c||}
 \hline
 Schl\"afli symbol & $Vol(\mathcal{B}_i \cap \hat{\mathcal{S}}(q,r))$ & $Vol(\hat{\mathcal{S}}(q,r))$ &  $\delta_{opt}(\mathcal{B}({q,r}))$.  \\ 
 \hline\hline
 $(\infty,3,3,\infty),~i=2$ & 0.1250000   & 0.1526609  &0.8188080  \\ 
 \hline
 $(\infty,3,4,\infty),~i=2$ & 0.1767766 & 0.2509603 & 0.7044011 \\
 \hline
 $(\infty,3,5,\infty),~i=2$ & 0.2022543 & 0.3323272 & 0.6085997 \\
 \hline
 $(\infty,3,6,\infty),~i=2$ & 0.2165064 & 0.4228923 & 0.5048035 \\
 \hline
 $(\infty,3,6,\infty),~i=0$ & 0.1443376 & 0.4228923 & 0.3365357 \\
 \hline
 $(\infty,4,3,\infty),~i=2$ & 0.1767766 & 0.2509603 & 0.7044011 \\
 \hline
 $(\infty,4,4,\infty),~i=2$ & 0.2500000 & 0.4579828 & 0.5458720 \\
 \hline
 $(\infty,4,4,\infty),~i=0$ & 0.2500000 & 0.4579828 & 0.5458720 \\
 \hline
 $(\infty,5,3,\infty),~i=2$ & 0.2022543 & 0.3323272 & 0.6085997 \\ 
 \hline
 $(\infty,6,3,\infty),~i=2$ & 0.2165064 & 0.4288923 & 0.5048035 \\
 \hline
 $(\infty,6,3,\infty),~i=0$ & 0.1443376 & 0.4288923 & 0.3365357 \\
 \hline
\end{tabular}
\end{center}
We summarize the results of this section in the following 
\begin{theorem}
In hyperbolic space $\mathbb{H}^3$, between the congruent horoball packings of one horoball, generated by simply truncated Coxeter orthoschemes with parallel faces, 
the $\mathcal{B}{(3,3)}$ horoball configuration provides the densest packing with density 
$\approx 0.8188080$.
\end{theorem}
\begin{figure}[h!]
  \begin{minipage}[b]{0.46\textwidth}
   \includegraphics[width=\textwidth]{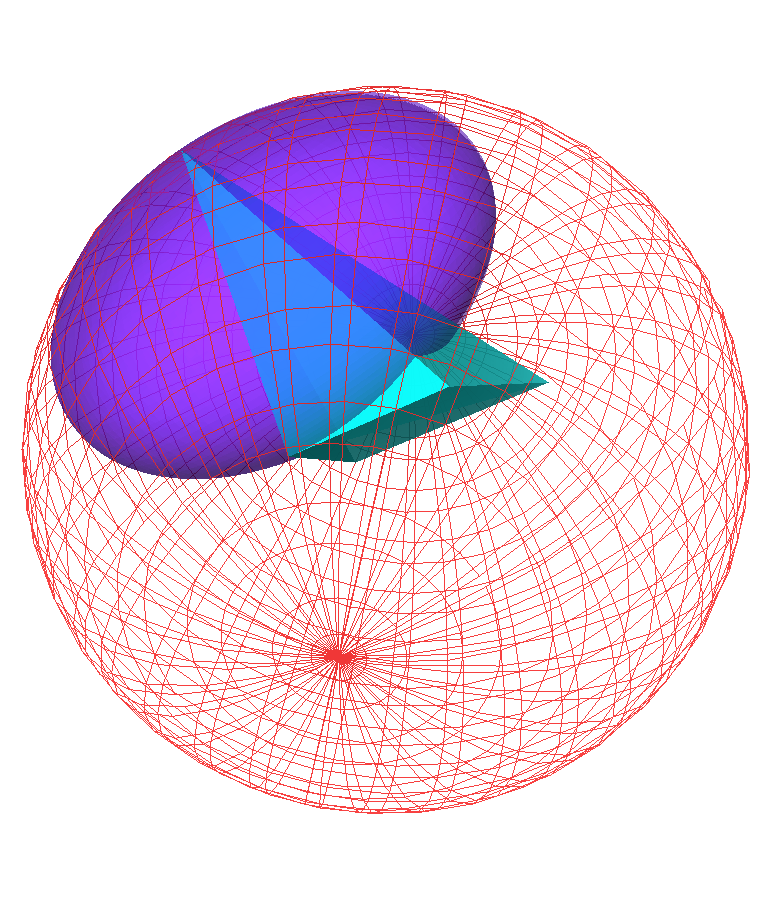}
   %\vspace{0.1 cm}
    \caption{The Horospheres related to the truncated orthoscheme $(\infty, 3, 3,\infty)$}
    \label{Horo33a}
    \end{minipage}
    \hspace{1 cm}
    \begin{minipage}[b]{0.46\textwidth}
    \includegraphics[width=\textwidth]{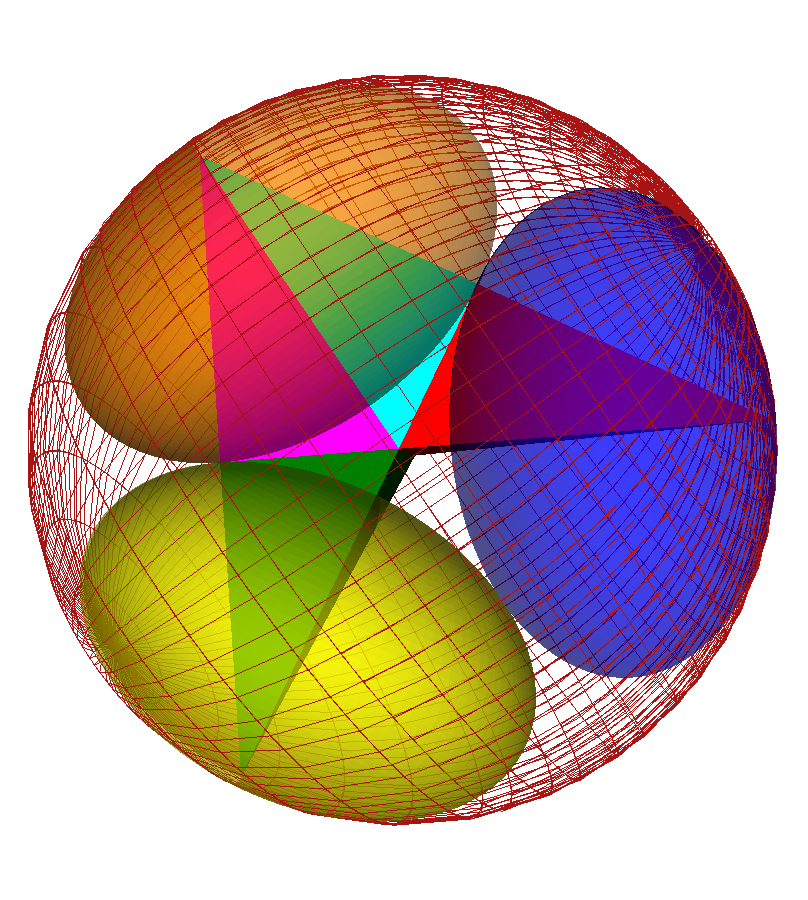}
    \caption{Horospheres packing configuration related to the truncated orthoscheme $\{ \infty, 3, 3,\infty \}$}
    \label{Horo33T}
  \end{minipage}
\end{figure}
\subsubsection{Horoball packings with two horoball types}
We are now focusing on the orthoschemes with the Schl\"afli symbols 
$\{ \infty, 3, 6,\infty \}$, 
$\{ \infty,4,4,\infty \}$, $\{ \infty,6,3,\infty \}$.
\begin{figure}[ht]
  \begin{minipage}[b]{0.48\textwidth}
    \includegraphics[width=\textwidth]{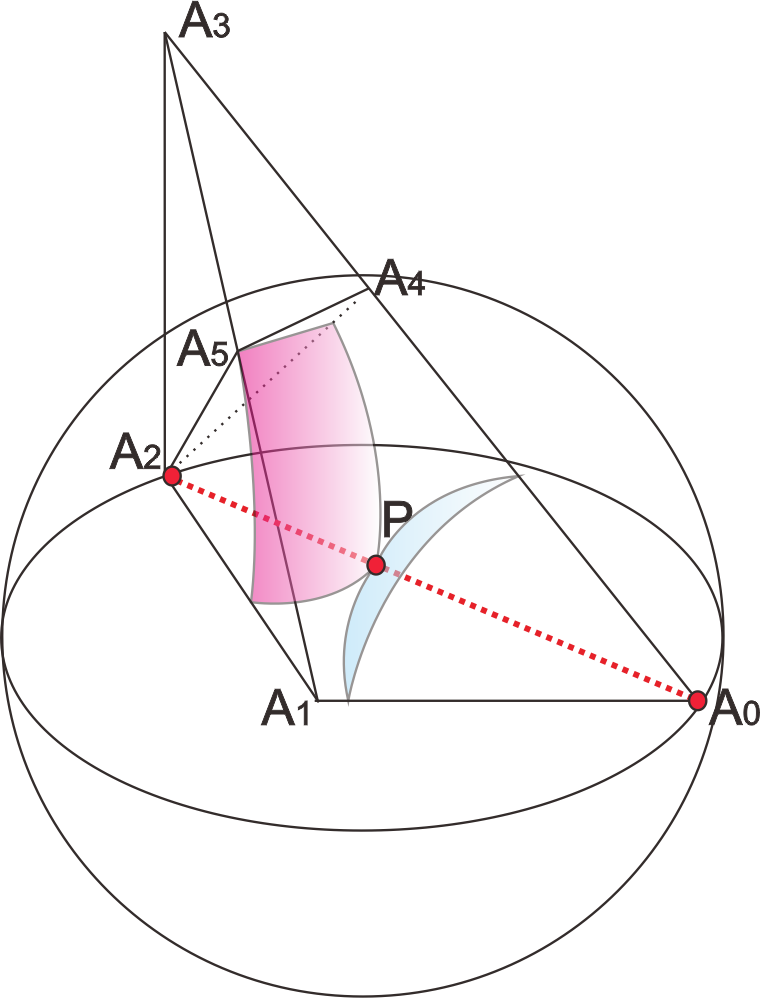}
    \vspace{0.1 cm}
    \caption{Sketch of two touching horospheres related to the truncated orthoscheme $\{\infty, 3, 6,\infty\}$}
    \label{doubleideal}
  \end{minipage}
  \hspace{0.5 cm}
  \begin{minipage}[b]{0.48\textwidth}
   \includegraphics[width=\textwidth]{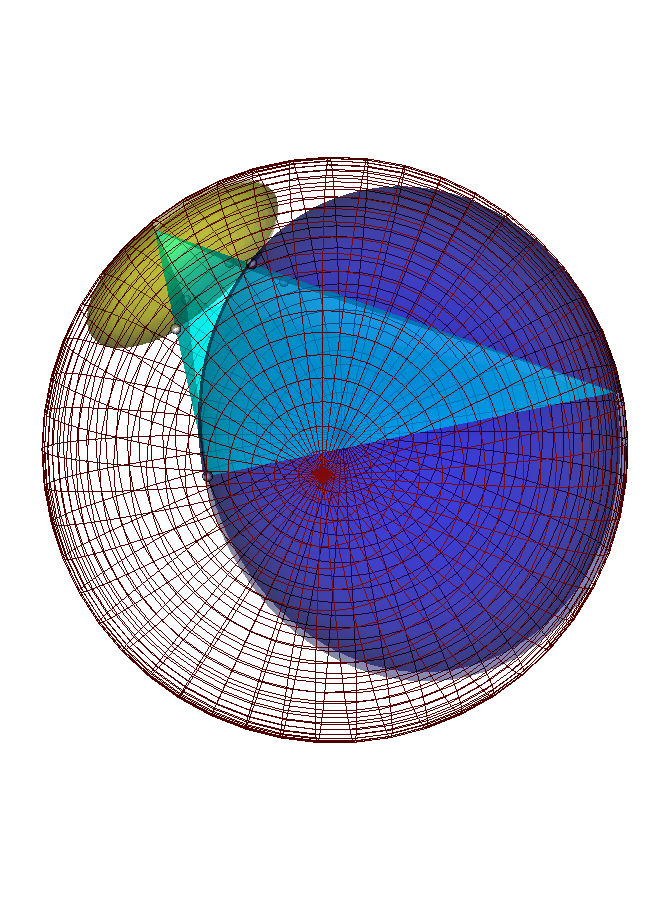}
    \caption{Two touching horospheres related to the truncated orthoscheme $\{\infty, 3, 6,\infty\}$}
    \label{Horo36}
    \end{minipage}
\end{figure}
In cases when the Coxeter simplex has multiple asymptotic vertices 
we allow horoballs of different types at different vertices. 

Two horoballs of a horoball packing are said to be of the same type or equipack-ed if and
only if their local packing densities with respect to a 
particular cell (in our case a truncated Coxeter
orthoscheme) are equal, otherwise the two horoballs are of different types. 
The set of all horoball types (they are congruent) at a vertex is a
one-parameter family. In our investigations we allow horoballs in different types 
(see \cite{KSz}, \cite{KSz2}, \cite{KSz3}).

As in the above subsection, first, we find the bounds for the largest possible
horoball type at each vertex lying at the infinity. Such a horoball is tangent to a face that not contains its center.
We set one horoball to be of the largest type, and increase the size of the other horoballs until they become tangent. 
We then vary the types of horoballs within the allowable bounds to find the optimal packing density. 
The following lemma proved e.g. in \cite{KSz2} gives a relationship between the volumes of two tangent horoball pieces 
centered at certain vertices of a tiling as we continuously vary their types.

Let $\tau_1$ and $\tau_2$ be two 
$n$-dimensional convex pyramid-like regions with vertices at $C_1$ and $C_2$ sharing a common edge $\overline{C_1C_2}$.
Let $B_1(x)$ and $B_2(x)$ denote two horoballs centered at $C_1$ and
$C_2$ tangent at the point
$I(x)\in {\overline{C_1C_2}}$. Define the point of tangency $I(0)$ (the ``midpoint") such that the
equality $V(0) = 2 Vol(B_1(0) \cap \tau_1) = 2 Vol(B_2(0) \cap \tau_2)$
holds for the volumes of the horoball sectors (see Fig.~6).

\begin{lemma}[\cite{KSz}]
\label{lemma:szirmai}
Let $x$ be the hyperbolic distance between $I(0)$ and $I(x)$,
then
\begin{equation}
V(x) = Vol(B_1(x) \cap \tau_1) + Vol(B_2(x) \cap \tau_2) = \frac{V(0)}{2}\left( e^{(n-1)x}+e^{-(n-1)x}\right) \notag
\end{equation}
strictly increases as $x\rightarrow\pm\infty$.
\end{lemma}

In the considered cases, the truncated orthoschemes have $2$ ideal vertices i.e $A_2$ and $A_0$ and we construct two horospheres in these vertices. 
The two horospheres touch each other at the point $P$ lying on the edge $A_0A_2$ (see Fig.~6). 
If we move the touching point along edge $A_0A_2$ the packing density is changing.
We can parameterize the possible movement of the touching point $P[\mathbf{p}]$ where $\mathbf{p}(t)=(1-t)\ba_2+t \ba_0$. Moreover, these two horospheres could
not intersect over their each opposite face, therefore there will be a restriction for parameter $t\in[t_1,t_2]$ (see Fig.~6). 
Then, for every possible $t$, we have to determine the parameters $s_i$ related to both horospheres $\mathcal{B}_i$ $(i=0,2)$.
\begin{enumerate}
    \item Constructing horosphere centered at $A_2$. 
    In our setting, the ideal point $A_2$ has the coordinate $(1,0,0,1)$, the horosphere could be immediately constructed 
    by applying formula (4.3). 
    To obtain its parameter $s_1$, we substitute the coordinates of point $P(t)$ to the equation (4.3) and then solve the expressed quadratic equation for $s_1$.
    \item Constructing a horosphere centred at $A_0$ is not straightforward because the horosphere equations in formula (4.3) holds only for horospheres centred at $A_2$. 
    However, we can use a transformation $T$ which is a hyperbolic isometry to transform the center $A_0$ and the to $A_2=(1,0,0,1)$ and we can apply 
    similar method as formerly substituting the coordinates 
    of $T(P(t))$ to the formula(4.3).
\end{enumerate}

The complete construction is presented in Fig.~6,~7.

Using the formula (4.1) we obtain the following formula for the density of a horoball packing $\mathcal{B}(q,r)$ related to Coxeter orthoscheme tiling $\cT(q,r)$:
\begin{equation}
\delta(\mathcal{B}({q,r,t}))=\frac{\sum_{i=1}^2 Vol(\mathcal{B}_i(t) \cap \hat{\mathcal{S}}(q,r))}{Vol(\hat{\mathcal{S}}(q,r))}
\end{equation}
where $\hat{\mathcal{S}}(q,r)$ denotes the fundamental domain of tiling $\cT(q,r)$,  
and $\mathcal{B}_i$ are the horoballs centered at the ideal vertices $A_0$ and $A_2$ $(q,r)=(3,6),(4,4),(6,3)$. 
Moreover, $t \in [t_1(q,r),t_2(q,r)]$ and this intervallum depends on parameters $q,r$ therefore the value of $\delta(\mathcal{B}({q,r,t}))$
is really depends on the parameters $t,q,r$.
\subsubsection{Horosphere packing density related to truncated orthoscheme\\ $\{\infty;3;6;\infty\}$ and $\{\infty;6;3;\infty\}$}
In this situations, there is in each case only one possible value of parameters $t$ $$t_{(3,6)}\approx 0.2119416,~ t_{(6,3)} \approx 0.5745582.$$ 
If $(q,r)=(3,6)$ then the horosphere $\cB_2$ touches the plane $u_2$ and $\cB_0$ touches the face $u_0$ and if 
$(q,r)=(6,3)$ $\cB_2$ touches the plane $u_2$ and $\cB_0$ touches the face $u_4$.

Finally, we obtain the volumes of horoball sectors and the optimal packing densities: 
\begin{equation*}
    \delta_{opt}(\mathcal{B}({3,6,t_{(3,6)}}))=\delta_{opt}(\mathcal{B}({6,3,t_{(6,3)}})) \approx 0.8413392.
\end{equation*}
 \subsubsection{Horosphere packing density related to truncated orthoscheme \\ $\{ \infty,4;4;\infty\} $}
 In this case, we obtain that the possible values of $t$ are $t\in [\approx 0.2150 < t < \approx 0.3497]$.
 We obtain the volumes of horoball sectors as the functions of $t$. 
 It is just similar to the previous case, the volume function of horoball sectors centered at $A_2$ is increasing function 
 of $t$ if the touching point moving with direction $A_2$ to $A_0$. 
 While the volume function of horoball sectors centered at $A_0$ is decreasing in this situation.
 
\begin{figure}[ht!]
\begin{center}
\includegraphics[scale=0.6]{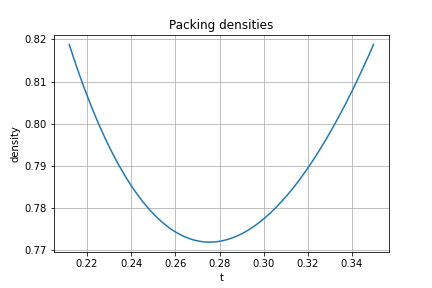}
\caption{The plot of density function for all possible $t$}
\label{dens44}
\end{center}
\end{figure}

In this case, we realize that the density is increasing as a function of $t$, see Fig.\ref{dens44}. Furthermore, the maximum density is attained when $t$ is largest, i.e when the horosphere centered $A_2$ touches the opposite face.
\begin{equation*}
\delta_{opt}(\mathcal{B}({4,4,t})=\sup \delta(\mathcal{B}({4,4,t})) \approx 0.8188081.
\end{equation*}
The data of packings  with two horoball types are summarized in the following table:
\begin{center}
 Table 4.
 
 \vspace{2mm}
 \begin{tabular}{||c c c c||}
 \hline
 Schl\"afli symbol & $Vol(\mathcal{B}_i \cap \hat{\mathcal{S}}(q,r))$ & $Vol(\hat{\mathcal{S}}(q,r))$ &  $\delta_{opt}(\mathcal{B}({q,r}))$.  \\ 
 \hline\hline
  $(\infty,3,6,\infty)$ & 0.3608439 & 0.4228923 & 0.8413392 \\
  \hline
  $(\infty,4,4,\infty)$ & 0.3750000 & 0.4579828 & 0.8188081 \\
  \hline
  $(\infty,6,3,\infty)$ & 0.3608439 & 0.4288923 & 0.8413392 \\
 \hline
\end{tabular}
\end{center}

Finally, we summarize our results of the sphere (inball) packings and the horosphere (horoball) packings in the following  
\begin{theorem}
In hyperbolic space $\mathbb{H}^3$, between the congruent ball and horoball packings of at most two horoball types, generated by simply truncated Coxeter orthoschemes with parallel faces, 
the $\mathcal{B}{(3,6)}$ and $\mathcal{B}{(6,3)}$ ball configuration provides the densest packing with density 
$\approx 0.8413392$.
\end{theorem}

\end{document}